\documentclass[12pt]{amsart}
\usepackage{graphicx}
\usepackage{hyperref}
\usepackage{setspace}
\usepackage{enumerate}
\usepackage{amsmath,amsfonts,amssymb,amscd}

\newtheorem{thm}{Theorem}[section]
\newtheorem{cor}[thm]{Corollary}
\newtheorem{lem}[thm]{Lemma}

\theoremstyle{definition}

\theoremstyle{remark}

\numberwithin{equation}{section}
\begin{document}

\title[]{Exponential sum estimate for systems \\ including linear polynomials}

\author{Shuntaro Yamagishi}
\address{Department of Mathematics \& Statistics \\
Queen's University \\
Kingston, ON\\  K7L 3N6 \\
Canada}
\email{sy46@queensu.ca}
\indent

\date{Revised on \today}

\begin{abstract}
In his paper \cite{S}, W. M. Schmidt obtained an exponential sum estimate for systems of polynomials
not including linear polynomials, which was then used to apply the Hardy-Littlewood circle method.
We prove an analogous estimate for systems including linear polynomials.
\end{abstract}

%{Primary 11D72, 11P32 ; Secondary}


\maketitle

\section{Introduction}
\label{sec intro}

\iffalse
Let us denote $\mathfrak{B}_0 = [0,1]^n$.
Let $\boldsymbol{\alpha} = (\boldsymbol{\alpha}_d, \ldots, \boldsymbol{\alpha}_1) \in \mathbb{R}^R$,
where $R = r_1 + \ldots + r_d$, and $\boldsymbol{\alpha}_{\ell} = (\alpha_{\ell, 1}, \ldots , \alpha_{\ell, r_{\ell}} ) \in \mathbb{R}^{r_{\ell}}$ $(1 \leq \ell \leq d)$.
We define
$$
\| \boldsymbol{\alpha} \| = \max_{ \substack{  1 \leq \ell \leq d  \\ 1 \leq r \leq r_{\ell}}} \| \alpha_{\ell, r} \|
\ \ \
\text{  and  } \ \  \
| \boldsymbol{\alpha} | = \max_{ \substack{  1 \leq \ell \leq d  \\ 1 \leq r \leq r_{\ell}}} | \alpha_{\ell, r} |.
$$
\fi

Let $\mathbf{u} = ( \mathbf{u}_{d}, \ldots, \mathbf{u}_{1} )$ be a system of polynomials in $\mathbb{Q}[x_1, \ldots, x_n]$, where
$\mathbf{u}_{\ell}  = ( u_{\ell,1}, \ldots, u_{\ell, r_{\ell}} )$ is the degree $\ell$ polynomials of $\mathbf{u}$  $(1 \leq \ell \leq d)$.
We let $\mathbf{U} = ( \mathbf{U}_{d}, \ldots, \mathbf{U}_{1} )$ be the system of forms, where for each $1 \leq \ell \leq d$,
$\mathbf{U}_{\ell} = ( U_{\ell,1}, \ldots, U_{\ell, r_{\ell}} )$
and $U_{\ell, r}$ is the degree $\ell$ portion of $u_{\ell, r}$ $(1 \leq r \leq r_{\ell} )$. Let us denote $\mathfrak{B}_0 = [0,1]^n$.
We define the following exponential sum associated to  $\mathbf{u}$,
\begin{equation}
\label{def of S 1}
S( \boldsymbol{\alpha}) = S( \mathbf{u},  \mathfrak{B}_0 ;\boldsymbol{\alpha}) := \sum_{\mathbf{x} \in P \mathfrak{B}_0 \cap \mathbb{Z}^n}
e \left( \sum_{1 \leq \ell \leq d} \sum_{ 1 \leq r \leq r_{\ell} } {\alpha}_{\ell, r} \cdot {u}_{\ell, r}  (\mathbf{x})  \right).
\end{equation}

In his paper \cite{S}, W. M. Schmidt obtained an exponential sum estimate for $S( \boldsymbol{\alpha})$
when $\mathbf{u}$ has integer coefficients, does not include linear polynomials, and satisfies certain properties.
The estimate was then used in applying the Hardy-Littlewood
circle method to obtain the asymptotic formula for the number of integer points of bounded height on the affine variety defined by
$\mathbf{u}$. We refer the reader to \cite{S} for more details on this important work.
The work of Schmidt was found useful in the breakthrough of B. Cook and {\'A}. Magyar \cite{CM}, where they count the number of solutions whose coordinates are all prime to diophantine equations, and also in \cite{XY}. It makes sense for Schmidt in \cite{S} to only consider
systems without linear polynomials, because he is concerned with integer points and linear polynomials can be eliminated via substitution in this case.
However, if one wants to apply the result of Schmidt for a coordinate dependent problem (where one can not eliminate linear polynomials by substitution),
then it may be useful to have analogous exponential sum estimates for systems including linear polynomials, and this is what we achieve in this paper.

We need to introduce some notations before we can state our result.
Let $1 < \ell \leq d$ and $r_{\ell} > 0$. We let $\mathbb{M}_{\ell} = \mathbb{M}_{\ell} (\mathbf{U}_{\ell})$ be
the affine variety in $(\mathbb{C}^n)^{\ell-1}$ associated to $\mathbf{U}_{\ell}$, for which the definition
we provide in (\ref{def of MU}) of Section \ref{appendix B}. For $R_0>0$, we denote $z_{R_0} (\mathbb{M}_{\ell})$ to be the number of integer points $(\mathbf{x}_1, \ldots, \mathbf{x}_{\ell-1} )$ on
$\mathbb{M}_{\ell}$ such that $$\max_{1 \leq i \leq \ell -1} \max_{1 \leq j \leq n}  | x_{ij} | \leq R_0,$$
where $\mathbf{x}_i = (x_{i1}, \ldots, x_{in}) \ (1 \leq i \leq \ell-1)$. We define $g_{\ell}( \mathbf{U}_{\ell} )$
to be the largest real number such that
\begin{equation}
\label{def gd+}
z_P(\mathbb{M}_{\ell}) \ll P^{n({\ell}-1) - g_{\ell}( \mathbf{U}_{\ell} ) + \varepsilon}
\end{equation}
holds for each $\varepsilon >0$.

%It was proved in \cite[pp. 280, Corollary]{S} that
%\begin{equation} \label{h and g+}
%h_{\ell}( \mathbf{U}_{\ell} ) < \frac{\ell!}{  (\log 2)^{\ell} }  \left( g_{\ell}( \mathbf{U}_{\ell}  ) + ({\ell}-1)r_{\ell} (r_{\ell} - 1)  \right).
%\end{equation}
Let
$$
\gamma_{\ell} = \frac{2^{{\ell}-1} ({\ell}-1) r_{\ell}}{ g_{\ell}( \mathbf{U}_{\ell} ) }
$$
when $r_{\ell} >0$ and $g_{\ell}( \mathbf{U}_{\ell} ) > 0$.
We let
$\gamma_{\ell} = 0$ if $r_{\ell} = 0$, and let
$\gamma_{\ell} = + \infty$ if $r_{\ell} > 0$ and $g_{\ell}( \mathbf{U}_{\ell} ) = 0$.

%For $\ell$ with $r_{\ell}>0$, we also define
%\begin{equation}
%\label{def gamma'+}
%\gamma'_{\ell} = \frac{ 2^{{\ell}-1} }{ g_{\ell}( \mathbf{U}_{\ell} ) } = \frac{ \gamma_{\ell} }{ ({\ell}-1) r_{\ell} }.
%\end{equation}

These quantities are not defined for linear polynomials. When $\ell = 1$, following \cite{CM} we define $\mathcal{B}_{1} (\mathbf{u}_1)$ to be the minimum number of non-zero coefficients in a non-trivial linear combination
$$
\lambda_1 U_1 + \ldots + \lambda_{r_1} U_{r_1},
$$
where $\boldsymbol{\lambda} = (\lambda_1, \ldots, \lambda_{r_1}) \in \mathbb{Q}^{r_1} \backslash \{ \mathbf{0} \}.$
Clearly $\mathcal{B}_{1} (\mathbf{u}_1) > 0$ if and only if the linear forms $U_1, \ldots, U_{r_1}$ are linearly independent over $\mathbb{Q}$.
If $r_1=0$ then we let $\mathcal{B}_{1} (\mathbf{u}_1) = + \infty$.

The following theorem is the main result of this paper.
\begin{thm}
\label{prop II in S}
Suppose $\mathbf{u}$ has coefficients in $\mathbb{Z}$, and that
$$
\mathcal{B}_1(\mathbf{u}_1) > 2 r_1 \left( \max \Big{ \{ } \
4 (r_1 + 1) \left( \sum_{j = 2}^d 4^{j- 2} \gamma_j \right)  , \  \frac{1}{4(R+1)} \
\Big{ \} }
\right)^{-1}.
$$
Let
$$
0 < \Omega < \min \Big{ \{ } \  \frac{1}{ 8r_1 + 9} \left( \sum_{j = 2}^d 4^{j- 2} \gamma_j \right)^{-1}, \  \left(  \frac{1}{2 (R + 1)} + \sum_{j = 2}^d 4^{j- 2} \gamma_j  \right)^{-1}  \ \Big{ \} }.
$$
Let $0 < \Delta \leq 1$, and let $P$ be sufficiently large with respect to $n$, $d$, $r_d, \ldots, r_1$, $\Delta$, $\Omega$, and $\mathbf{u}$.
Then one of the following two alternatives must hold:

$(i)$ $|S(\boldsymbol{\alpha})| \leq P^{n - \Delta \Omega}$.

$(ii)$ There exists $q \in \mathbb{N}$ such that
$$
q \leq P^{\Delta} \ \ \text{ and }  \  \  \|  q \boldsymbol{\alpha}_{\ell}  \| \leq P^{ - \ell + \Delta} \ \ (1 \leq \ell \leq d).
$$
%\newline The implicit constants depend at most on $n,d, r_d, \ldots, r_1, \eta, \varepsilon$, and $\mathbf{u}$.
\end{thm}

In Section \ref{appendix B}, we also prove a lemma on estimating the quantity known as the singular integral, which comes up
%appears 
in the Hardy-Littlewood circle method. We use $\ll$ and $\gg$ to denote Vinogradov's well-known notation.
We also use the notation $e(x)$ to denote $e^{2\pi i x}$.
For $q \in \mathbb{N}$, we use the numbers from $\{0, 1, \ldots, q-1 \}$ to represent the residue classes of
$\mathbb{Z}/q\mathbb{Z}$.

\section{Proof of Theorem \ref{prop II in S}}
\label{appendix B}

First we present the following lemma from \cite{CM}.
\begin{lem} \cite[Lemma 3]{CM}
\label{Lemma 3 in CM}
Let $\mathbf{G} = (G_1, \ldots, G_{r'})$ be a system of linear forms in $\mathbb{Q}[x_1, \ldots, x_n]$.
Given any $1 \leq j \leq n$, we have
$$
\mathcal{B}_1(\mathbf{G}|_{x_j = 0}) \geq \mathcal{B}_1(\mathbf{G}) - 1.
$$
\end{lem}

%Let us denote %$\mathfrak{B}_1 = [-1,1]^n$ and
%$\mathfrak{B}_0 = [0,1]^n$.
Let $\mathbf{x} = (x_1, \ldots, x_n)$ and $\mathbf{x}_j = (x_{j,1}, \ldots, x_{j,n})$ for $j \geq 1$.
Given a function $G(\mathbf{x})$, we define
$$
\Gamma_{\ell, G} (\mathbf{x}_1, \ldots, \mathbf{x}_{\ell}) = \sum_{t_1=0}^1 \ldots \sum_{t_{\ell}=0}^1 (-1)^{t_1 + \ldots + t_{\ell}} \
G( t_1 \mathbf{x}_1 + \ldots+ t_{\ell} \mathbf{x}_{\ell} ).
$$
Then it follows that $\Gamma_{{\ell}, G}$ is symmetric in its ${\ell}$ arguments, and that
$\Gamma_{{\ell}, G} (\mathbf{x}_1, \ldots,\mathbf{x}_{{\ell}-1}, \mathbf{0}) = 0$ \cite[Section 11]{S}.
%It is clear from the definition that if $G'(\mathbf{x})$ is another function, then $\Gamma_{{\ell}, G} + \Gamma_{{\ell}, G'} = \Gamma_{{\ell}, G + G'}.$
We also have that if $G$ is a form of degree $d$ and ${\ell} > d > 0$, then $\Gamma_{{\ell}, G}= 0$ \cite[Lemma 11.2]{S}.

For $\alpha \in \mathbb{R}$, let $\| \alpha \|$ denote the
distance from $\alpha$ to the closest integer. Let $\boldsymbol{\alpha} = (\boldsymbol{\alpha}_d, \ldots, \boldsymbol{\alpha}_1) \in \mathbb{R}^R$,
where $R = r_1 + \ldots + r_d$ and $\boldsymbol{\alpha}_{\ell} = (\alpha_{\ell, 1}, \ldots , \alpha_{\ell, r_{\ell}} ) \in \mathbb{R}^{r_{\ell}}$ $(1 \leq \ell \leq d)$.
We define
$$
\| \boldsymbol{\alpha} \| = \max_{ \substack{  1 \leq \ell \leq d  \\ 1 \leq r \leq r_{\ell}}} \| \alpha_{\ell, r} \|
\ \ \ \text{ and } \ \ \
| \boldsymbol{\alpha} | = \max_{ \substack{  1 \leq \ell \leq d  \\ 1 \leq r \leq r_{\ell}}} | \alpha_{\ell, r} |.
$$

Let $\mathbf{u} = ( \mathbf{u}_{d}, \ldots, \mathbf{u}_{1} )$ and  $\mathbf{U} = ( \mathbf{U}_{d}, \ldots, \mathbf{U}_{1} )$ be
as in Section \ref{sec intro}.
Let $\mathbf{e}_1, \ldots, \mathbf{e}_n$ be the standard basis vectors of $\mathbb{C}^n$. Let $1 < \ell \leq d$.
We define $\mathbb{M}_{\ell} = \mathbb{M}_{\ell} (\mathbf{U}_{\ell}) $ to be the set of $(\ell-1)$-tuples $(\mathbf{x}_1, \ldots, \mathbf{x}_{\ell-1} ) \in (\mathbb{C}^n)^{\ell-1}$ for which the matrix
\begin{equation}
\label{def of MU}
[m_{i r}] = [ \Gamma_{\ell, U_{\ell, r}} ( \mathbf{x}_1, \ldots , \mathbf{x}_{d-1}, \mathbf{e}_i ) ] \ \ \ \  (1 \leq r \leq r_{\ell}, 1\leq i \leq n)
\end{equation}
has rank strictly less than $r_{\ell}$.

Lemma \ref{Lemma 15.1 in S} below is the inhomogeneous polynomials version of \cite[Lemma 15.1]{S},
and it is obtained by essentially the same proof. We refer the reader to \cite[Section 9]{S} and `Remark on inhomogeneous polynomials' in \cite[pp. 262]{S}
for further explanation. We remark that the implicit constants may depend on $\mathbf{u}$ here, and not only on $\mathbf{U}$.
We also note that \cite[Lemma 15.1]{S} is for systems without linear polynomials in contrast to Lemma \ref{Lemma 15.1 in S} below.
However, it is clear from the proof of \cite[Lemma 15.1]{S} that the lemma is not affected with the presence of linear polynomials.
%If we have Condition $(\star)$, where the coefficients of lower order of $u_{\ell, r} (\mathbf{x})$ depend on $P$
%in a controlled manner, then we can obtain more than just for the polynomials with the highest degree in the system.
\begin{lem}\cite[Lemma 15.1]{S}
\label{Lemma 15.1 in S}
Suppose $\mathbf{u}$ has coefficients in $\mathbb{Z}$.
Let $Q > 0$ and  $\varepsilon >0$. Let $2 \leq \ell \leq d$ with $r_{\ell} > 0$. Let $P$ be sufficiently large with respect
to $d$ and  $r_d, \ldots, r_1$.
If $\ell =d$, then let $\theta = 0$ and $q=1$.
On the other hand, if
$2 \leq \ell < d$, then suppose $0 \leq \theta < 1/4$ and that there is $q \in \mathbb{N}$ with
$$
q \leq P^{\theta} \ \ \text{ and } \ \  \| q \boldsymbol{\alpha}_{j}  \| \leq P^{\theta - j} \ \ (\ell < j \leq d).
$$
Let $S( \boldsymbol{\alpha})$ be the sum associated to $\mathbf{u}$ as in ~(\ref{def of S 1}). %Given $\eta > 0$ such that $\eta + 4 \theta \leq 1$,
Given $\eta > 0$ with $\eta + 4 \theta \leq 1$, one of the following three alternatives must hold:

$(i)$ $|  S( \boldsymbol{\alpha})  | \leq P^{n-Q}$.

$(ii)$ There exists $n_0 \in \mathbb{N}$ such that
$$
n_0 \ll P^{r_{\ell}(\ell-1) \eta} \text{  and  } \|  q n_0  \boldsymbol{\alpha}_{\ell} \| \ll P^{ -\ell + 4 \theta + r_{\ell} (\ell-1) \eta}.
$$

$(iii)$ $ z_{R_0} (\mathbb{M}_{\ell}) \gg R_0^{ (\ell-1)n - 2^{\ell-1}(Q/ \eta) - \varepsilon}$
holds with $R_0 = P^{\eta}$.
\newline
\newline
%When $\ell = d$, the implicit constants depend at most on $n,d, r_d, \eta, \varepsilon$ and $\mathbf{U}^{(d)}$.
%Otherwise, when $2 \leq \ell < d$,
The implicit constants
depend at most on $n,d, r_d, \ldots, r_1, \eta, \varepsilon$, and $\mathbf{u}$.
%and all the $D^{(j)}_{\ell, r}$ from Condition $(\star)$.
\end{lem}

We are left to deal with the case $\ell = 1$ in Lemma \ref{Lemma 15.1 in S}.
Given $\boldsymbol{\epsilon} \in (\mathbb{N} \cup \{0\} )^n$ and sufficiently differentiable function $f: \mathbb{R}^n \rightarrow
\mathbb{C}$, put
$$
\partial^{\boldsymbol{\epsilon} } f = \frac{\partial^{\epsilon_1 + \ldots + \epsilon_n} f }{ \partial x_1^{\epsilon_1} \ldots \  \partial x_n^{\epsilon_n}}.
$$
Let $\mathcal{C}^n(\mathbb{R}^n)$ be the set of $n$-th continuously differentiable functions defined on $\mathbb{R}^n$.

For $\boldsymbol{\epsilon} \in \{ 0, 1\}^n$, we define $\overline{\boldsymbol{\epsilon}} = (1, 1, \ldots , 1) - \boldsymbol{\epsilon}$.
Given $\mathbf{t} = (t_1, \ldots, t_n)$, we let $\mathbf{t}_{\boldsymbol{\epsilon}}$ be the vector whose $i$-th coordinate equals zero if $\epsilon_i = 0$ and equals $t_i$ if $\epsilon_i = 1$. Similarly, given $\mathbf{N} = (N_1, \ldots, N_n) \in \mathbb{Z}^n$, we let
$\mathbf{N}_{ \overline{\boldsymbol{\epsilon}} }$
be the vector whose $i$-th coordinate equals $N_i$ if $\epsilon_i = 0$ and equals zero if $\epsilon_i = 1$.
The following is a generalization of the partial summation formula obtained by applying induction on the dimension.
\begin{lem}\cite[Lemma 2.1]{BP}
\label{partial sum}
Let $\varrho : \mathbb{Z}^n \rightarrow \mathbb{C}$ be a function, and
let
$$
T_{\varrho}(\mathbf{t}) = \sum_{0 \leq x_1 \leq t_1} \ldots \sum_{0 \leq x_n \leq t_n} \varrho(\mathbf{x}).
$$
Then for any $f \in \mathcal{C}^n(\mathbb{R}^n)$ we have
\begin{eqnarray}
\sum_{  \substack{  0 \leq x_i \leq N_i \\ (1 \leq i \leq n) }  } f(\mathbf{x}) \varrho(\mathbf{x})
&=& \sum_{ \boldsymbol{\epsilon} \in \{ 0, 1\}^n  }
\left( \prod_{1 \leq i \leq n} (-1)^{\epsilon_i} N_i^{\epsilon_i - 1} \right) \cdot
\\
&& \int_{[0,N_1]} \ldots  \int_{[0,N_n]} \partial^{\boldsymbol{\epsilon} } f ( \mathbf{N}_{ \overline{\boldsymbol{\epsilon}} } + \mathbf{t}_{\boldsymbol{\epsilon}} )
\ T_{\varrho}( \mathbf{N}_{ \overline{\boldsymbol{\epsilon}} } + \mathbf{t}_{\boldsymbol{\epsilon}} ) \ d t_n \ldots d t_1.
\notag
\end{eqnarray}

\end{lem}

Let us use the following notations. For $\mathbf{a} = (a_1, \ldots, a_{r_1}) \in (\mathbb{Z}/ q \mathbb{Z})^{r_1}$, we let
$$
\mathfrak{M}_{\mathbf{a}, q}(C) = \{ \boldsymbol{\alpha}_1 \in  [0,1)^{r_1} :  \max_{1 \leq r \leq r_1}| \alpha_{1,r} - a_r / q| \leq P^{C - 1}    \},
$$
$$
\mathfrak{M}(C) = \bigcup_{\substack{ \gcd(\mathbf{a}, q) = 1 \\ \mathbf{a} \in (\mathbb{Z}/ q \mathbb{Z})^{r_1} \\ 1 \leq q \leq P^C }} \mathfrak{M}_{\mathbf{a}, q}(C),
$$
and
$$
\mathfrak{m}(C) =  [0,1)^{r_1}  \backslash \mathfrak{M}(C).
$$
We also let
$$
\mathfrak{N}_{a, q}(C) = \{ {\alpha} \in  [0,1):  | \alpha - a / q| \leq P^{C - 1}    \},
$$
$$
\mathfrak{N}(C) = \bigcup_{\substack{ \gcd(a, q) = 1  \\  0 \leq a < q \\ 1 \leq q \leq P^C }} \mathfrak{N}_{a, q}(C),
$$
and
$$
\mathfrak{n}(C) =  [0,1)  \backslash \mathfrak{N}(C).
$$

With the use of Lemma \ref{partial sum}, we obtain the following result when $r_1 > 0$.
\begin{lem}
\label{Lemma 15.1 in S-linear}
Suppose $\mathbf{u}$ has coefficients in $\mathbb{Z}$ and that $r_1 > 0$. Let $0 < \theta_0 < 1$, and suppose there exists $q \in \mathbb{N}$ with
$$
q \leq P^{\theta_0} \ \ \text{ and } \ \  \| q \boldsymbol{\alpha}_{j}  \| \leq P^{\theta_0 - j} \ \ (1 < j \leq d).
$$
Let $S( \boldsymbol{\alpha})$ be the sum associated to $\mathbf{u}$ as in ~(\ref{def of S 1}). %Given $\eta > 0$ such that $\eta + 4 \theta \leq 1$,
Let $\varepsilon_0 > 0$ be sufficiently small. Let  $Q > 0$ and $0<  Q_0 < 1$ be two real numbers such that
$$
\theta_0 < \frac{Q_0/2 - \varepsilon_0}{2 r_1}
$$
and
$$
Q <   \mathcal{B}_1(\mathbf{u}_1) \left(   \frac{Q_0/2 - \varepsilon_0}{r_1}  - 2 \theta_0 \right).
$$
Suppose $P$ is sufficiently large with respect to $d$, $n$, $r_d, \ldots, r_1$, $\varepsilon_0$, $\theta_0$, $Q_0$, $Q$, and $\mathbf{u}$. Then one of the following two alternatives must hold:

$(i)$ $|  S( \boldsymbol{\alpha})  | \leq P^{n-Q}$.

$(ii)$ There exists $n_0 \in \mathbb{N}$ such that
$$
n_0 \leq P^{Q_0} \text{  and  } \|  n_0  \boldsymbol{\alpha}_{1} \| \leq P^{ Q_0 - 1 }.
$$
\end{lem}
\begin{proof}
If the alternative $(ii)$ holds then we are done. Thus let us assume it is not the case.
Suppose $\boldsymbol{\alpha}_1 \in \mathfrak{M}(Q_0/2)$. Then for some
$1 \leq q' \leq P^{Q_0/2}$ and $a_1, \ldots, a_{r_1} \in \mathbb{Z}$, we have
$$
\max_{1 \leq r \leq r_1} | \alpha_{1,r} - a_{r}/q' | \leq P^{(Q_0/2) - 1}
$$
from which it follows that
$$
\| q' \boldsymbol{\alpha}_{1}  \| \leq P^{Q_0- 1},
$$
and this is a contradiction. Therefore, we have $\boldsymbol{\alpha}_1 \in \mathfrak{m}(Q_0/2)$.

For simplicity we denote $B= \mathcal{B}_1(\mathbf{u}_1)$ and $Q'_0 = Q_0/2$.
Let us also denote
$$
\sum_{r=1}^{r_1} {\alpha}_{1, r} \cdot {U}_{1, r}  (\mathbf{x}) =  \gamma_1  x_1 + \ldots +  \gamma_n  x_n.
$$
We let $\widetilde{M}_1$ be the $n \times r_1$ matrix, where its $(j,r)$-th entry is the $x_j$ coefficient of $U_{1,r}(\mathbf{x})$.
Since this matrix has full rank (because $B>0$), let us take an invertible $r_1 \times r_1$ minor, which we assume without loss of generality to be the first $r_1$ rows of $\widetilde{M}_1$, and denote it $M_1$.

Suppose $\gamma_1, \ldots, \gamma_{r_1} \in \mathfrak{N}^{}(C')$ for some $C'>0$. Then there exist integers $a_1, \ldots, a_{r_1}$ and $q_{1}, \ldots, q_{r_1}$ such that $\gcd({a}_r, q_r) = 1$, $0< q_r \leq P^{C'}$, and $| \gamma_{r} - a_{r} /q_r | \leq P^{C'}/ P$ $(1 \leq r \leq r_1)$.
Let us define
$$
\left[ {\begin{array}{c}
a'_{1} /q'\\ \vdots  \\  a'_{r_{1} } /q'
\end{array}} \right] =
M_1^{-1} \cdot \left[ {\begin{array}{c}
a_{1} /q_1 \\ \vdots  \\  a_{r_{1}} /q_{r_1}
\end{array}} \right]
\ \ \text{ and } \ \
\left[ {\begin{array}{c}
\beta'_{1} \\ \vdots  \\  \beta'_{r_{1}}
\end{array}} \right] =
M_1^{-1} \cdot \left[ {\begin{array}{c}
\gamma_{1} - a_{1}/ q_1 \\ \vdots  \\  \gamma_{r_{1}} - a_{r_{1}}/q_{r_1}
\end{array}} \right].
$$
It is easy to deduce that we have
$$
q' \leq P^{ r_1 C' + \varepsilon_0} \ \ \text{ and } \ \ |\beta'_{r}| \leq  \frac{P^{r_1 C' + \varepsilon_0}}{P} \ \ (1 \leq r \leq r_{1})
$$
when $P$ is sufficiently large with respect to the coefficients of $\mathbf{U}_1$.
Since $\alpha_{1,r} = \frac{a'_{r} }{q'} + \beta'_{r}$ $(1 \leq r \leq r_1)$, we see that $ \boldsymbol{\alpha}_1 \in \mathfrak{M}(r_1 C' + \varepsilon_0)$.
However, since $\boldsymbol{\alpha}_1 \in \mathfrak{m}(Q'_0)$, it follows from this argument that at least one of  $\gamma_1, \ldots, \gamma_{r_1}$ is in
$\mathfrak{n}( (Q'_0 -\varepsilon_0)/r_1)$. Without loss of generality, we suppose that $\gamma_1 \in \mathfrak{n}( (Q'_0 -\varepsilon_0)/r_1)$.

Let $\widetilde{M}_2$ be the matrix obtained by removing the first row of $\widetilde{M}_1$.
If $B-1> 0$, then we know that $\widetilde{M}_2$ has full rank. Let us take an invertible $r_1 \times r_1$ minor, which we assume without loss of generality to be the first $r_1$ rows of $\widetilde{M}_2$, and denote it $M_2$. By the same argument as above, we obtain without loss of generality that $\gamma_2 \in \mathfrak{n}( (Q'_0-\varepsilon_0)/r_1)$.
In fact we can repeat the argument $B$ times, and obtain
that $\gamma_1, \gamma_2, \ldots,  \gamma_B \in \mathfrak{n}( (Q'_0 -\varepsilon_0)/r_1)$.

Since
$q \leq P^{\theta_0} \leq P^{(Q'_0-\varepsilon_0)/r_1}$, it then follows that
\begin{equation}
\label{lower bound with theta0}
\frac{P^{(Q'_0-\varepsilon_0)/r_1}}{P} <  \| q \gamma_i \|  \ \ (1 \leq i \leq B).
\end{equation}

For each $2 \leq \ell \leq d, 1 \leq r \leq r_{\ell}$, let $a_{\ell, r} \in \mathbb{Z}$ and $\beta_{\ell, r} \in \mathbb{R}$ be such that
\begin{equation}
\label{bunch of alphas}
 \alpha_{\ell, r} - a_{\ell, r} / q = \beta_{\ell, r} \  \  \text{ and } \  \  | \beta_{\ell, r}| \leq P^{\theta_0 - \ell}.
\end{equation}

We then consider
\begin{eqnarray}
\label{bound in S with linear 0}
| S( \boldsymbol{\alpha}) | &=&  \Big{|} \sum_{ \substack{ 0 \leq k_i < q \\ (1 \leq i \leq n)} } \  \sum_{  \substack{ \mathbf{x} \in [0,P]^n
\\ x_i \equiv k_i (\text{mod }q) \\ (1 \leq i \leq n) } }
e \left( \sum_{1 \leq \ell \leq d} \sum_{ 1 \leq r \leq r_{\ell} } {\alpha}_{\ell, r} \cdot {u}_{\ell, r}  (\mathbf{x})  \right) \Big{|}
\\
&\leq& q^n  \max_{ \substack{ 0 \leq k_i < q \\ (1 \leq i \leq n)} }  \Big{|}  \sum_{ \substack{ 0 \leq y_i \leq (P - k_i)/ q  \\  (1 \leq i \leq n)  } }
e \left( \sum_{1 \leq \ell \leq d} \sum_{ 1 \leq r \leq r_{\ell} } {\alpha}_{\ell, r} \cdot {u}_{\ell, r}  (q \mathbf{y} + \mathbf{k} )  \right) \Big{|}.
\notag
\end{eqnarray}
Let us denote
$$
f(\mathbf{y}) = e \left( \sum_{2 \leq \ell \leq d} \sum_{ 1 \leq r \leq r_{\ell} } {\beta}_{\ell, r} \cdot {u}_{\ell, r}  (q \mathbf{y} + \mathbf{k} ) \right).
$$
Using the fact that $e(m) = 1$ for $m \in \mathbb{Z}$, we can simplify the above inequality (\ref{bound in S with linear 0}) further,
\begin{eqnarray}
\notag
| S( \boldsymbol{\alpha}) |
&\leq& q^n  \max_{ \substack{ 0 \leq k_i < q \\ (1 \leq i \leq n)} }  \Big{|}  \sum_{ \substack{ 0 \leq y_i \leq (P - k_i)/ q  \\  (1 \leq i \leq n)  } }
e \left( \sum_{2 \leq \ell \leq d} \sum_{ 1 \leq r \leq r_{\ell} } {\beta}_{\ell, r} \cdot {u}_{\ell, r}  (q \mathbf{y} + \mathbf{k} )
+ \sum_{ 1 \leq r \leq r_{1} } {\alpha}_{1, r} \cdot {U}_{1, r}  (q \mathbf{y} )
\right) \Big{|}
\\
&\leq& q^n  \max_{ \substack{ 0 \leq k_i < q \\ (1 \leq i \leq n)} } \ \sum_{ \substack{ 0 \leq y_i \leq (P - k_i)/ q  \\  (B < i \leq n)  } }
\Big{|}  \sum_{ \substack{ 0 \leq y_i \leq (P - k_i)/ q  \\  (1 \leq i \leq B)  } }
f(\mathbf{y}) \  e \left( \sum_{ 1 \leq i \leq B } q \gamma_i y_i
\right) \Big{|}.
\notag
\end{eqnarray}

Let $0 \leq y_i \leq (P - k_i)/ q $ $(B < i \leq n)$.
Given $\boldsymbol{\epsilon} \in \{ 0,1\}^B$, let
$\left( \frac{(P - \mathbf{k})}{q}  \right)_{ \overline{\boldsymbol{\epsilon}} }$
be the vector whose $i$-th coordinate, for $1 \leq i \leq B$, equals $(P - k_i)/q$ if $\epsilon_i = 0$ and equals zero if $\epsilon_i = 1$,
and for $B < i \leq n$, equals $y_i$. We also let $\mathbf{t}_{ \boldsymbol{\epsilon} }$
be the vector whose $i$-th coordinate, for $1 \leq i \leq B$, equals $0$ if $\epsilon_i = 0$ and equals $t_i$ if $\epsilon_i = 1$,
and for $B < i \leq n$, equals zero.

We prove that given $\boldsymbol{\epsilon} \in \{ 0,1\}^B$ and $0 \leq t_i \leq (P - k_i)/q$ $(1 \leq i \leq B)$, we have
\begin{equation}
\label{partial bound1}
\frac{\partial^{\epsilon_1 + \ldots + \epsilon_B} f }{ \partial y_1^{\epsilon_1} \ldots \  \partial y_B^{\epsilon_B}}
\Big{|}_{\mathbf{y} =   \left( \frac{(P - \mathbf{k})}{q}  \right)_{ \overline{\boldsymbol{\epsilon}} }  + \mathbf{t}_{\boldsymbol{\epsilon}} }
\ll
q^{\epsilon_1 + \ldots + \epsilon_B}  P^{ (\theta_0 - 1) (\epsilon_1 + \ldots + \epsilon_B) },
\end{equation}
where the implicit constant is independent of $k_1, \ldots, k_n$, $y_{B+1}, \ldots, y_n$, and $\mathbf{t}$.
In order to prove this statement, without loss of generality suppose $\epsilon_i = 1$ for $1 \leq i \leq E$
and $\epsilon_i = 0$ for $E < i \leq B$. The statement is trivial if $\epsilon_i = 0$ for all $1 \leq i \leq B$.
Let $i_1 < \ldots < i_m \leq E$.
First note when $m \leq d$, we have from (\ref{bunch of alphas}) that
\begin{eqnarray}
\label{partial bound1+}
&&\frac{\partial^{ m } }{ \partial y_{i_1}  \ldots \  \partial y_{i_m}  }
\left( \sum_{2 \leq \ell \leq d} \sum_{ 1 \leq r \leq r_{\ell} } {\beta}_{\ell, r} \cdot {u}_{\ell, r}  (q \mathbf{y} + \mathbf{k} ) \right)
\Big{|}_{\mathbf{y} =   \left( \frac{(P - \mathbf{k})}{q}  \right)_{ \overline{\boldsymbol{\epsilon}} }  + \mathbf{t}_{\boldsymbol{\epsilon}} }
\\
&\ll&
q^{ m  } \sum_{\max\{2,  m\} \leq \ell \leq d} \ \sum_{ 1 \leq r \leq r_{\ell} } {\beta}_{\ell, r} P^{\ell - m }
\notag
\\
&\ll&
q^m P^{ \theta_0 - m},
\notag
\end{eqnarray}
and when $m > d$,
\begin{eqnarray}
\label{partial bound1++}
\frac{\partial^{ m } }{ \partial y_{i_1}  \ldots \  \partial y_{i_m}  }
\left( \sum_{2 \leq \ell \leq d} \sum_{ 1 \leq r \leq r_{\ell} } {\beta}_{\ell, r} \cdot {u}_{\ell, r}  (q \mathbf{y} + \mathbf{k} ) \right) =   0.
\end{eqnarray}
Thus we have
\begin{eqnarray}
\label{partial bound1+++}
\frac{\partial^{ E } f}{ \partial y_{1}  \ldots \  \partial y_{E}  }  \Big{|}_{\mathbf{y} =   \left( \frac{(P - \mathbf{k})}{q}  \right)_{ \overline{\boldsymbol{\epsilon}} }  + \mathbf{t}_{\boldsymbol{\epsilon}}  }
&\ll&
\max_{ \substack{ m_1 + \ldots + m_j = E \\ 1 \leq m_i \leq d  \\ (1 \leq i \leq j) } } q^E P^{ j \theta_0 - E}
\\
&\ll&
q^E P^{ (\theta_0 - 1) E},
\notag
\end{eqnarray}
from which we can deduce (\ref{partial bound1}).
We now prepare to apply Lemma \ref{partial sum}.
Let $0 \leq t_i \leq (P - k_i) /q$ $(1 \leq i \leq B)$. It follows from (\ref{lower bound with theta0}) that
\begin{eqnarray}
\label{inequality of exp sum}
\Big{|} \sum_{ 0 \leq y_i \leq t_i } e( q \gamma_i  y_i ) \Big{|} \ll \min \Big{\{} \  t_i + 1, \  \| q \gamma_i \|^{-1}  \Big{\}}  \leq P^{1 - (Q'_0 - \varepsilon_0)/r_1} \ \ \ (1 \leq i \leq B).
\end{eqnarray}
Then for $\boldsymbol{\epsilon} \in \{0,1 \}^B,$ we have by (\ref{partial bound1}) and (\ref{inequality of exp sum}) that
%We then have by (\ref{partial bound1}),
\begin{eqnarray}
\label{partial integral}
&&\int_{[0,(P-k_1)/q]} \ldots  \int_{[0,(P - k_B)/q]} \partial^{\boldsymbol{\epsilon} } f \left( \left( \frac{(P - \mathbf{k})}{q}  \right)_{ \overline{\boldsymbol{\epsilon}} } +  \mathbf{t}_{\boldsymbol{\epsilon}} \right) \cdot
\\
&&\phantom{1234567890} \sum_{ \substack{ 0 \leq y_i \leq (P - k_i)/ q  \\  \epsilon_i = 0  } } \
\sum_{ \substack{ 0 \leq y_i \leq t_i \\  \epsilon_i = 1  } }
e \left(\sum_{ 1 \leq i \leq B } q \gamma_i y_i
\right) \ d t_B \ldots d t_1
\notag
\\
&\ll& q^{\epsilon_1 + \ldots + \epsilon_B} P^{(\theta_0 - 1) (\epsilon_1 + \ldots + \epsilon_B) } \left( \prod_{1 \leq i \leq B} \frac{ P- k_i}{q} \right) \cdot
P^{B - B(Q'_0 - \varepsilon_0)/r_1}.
\notag
\end{eqnarray}

Therefore, by Lemma \ref{partial sum} and (\ref{partial integral}) we obtain for any $0 \leq y_i \leq (P - k_i)/q$ $(B < i \leq n)$,
\begin{eqnarray}
\label{bound in S with linear 1}
&&\Big{|} \sum_{ \substack{ 0 \leq y_i \leq (P - k_i)/ q  \\  (1 \leq i \leq B)  } }
f(\mathbf{y}) \  e \left( \sum_{ 1 \leq i \leq B } q {\gamma}_{i} y_i
\right) \Big{|}
\\
&\ll&
\sum_{ \boldsymbol{\epsilon} \in \{ 0,1\}^B }
\left( \prod_{1 \leq i \leq B} \left(\frac{ P - k_i}{q} \right)^{\epsilon_i - 1} \right)
q^{\epsilon_1 + \ldots + \epsilon_B}  \cdot
\notag
\\
\notag
&& P^{(\theta_0 - 1) (\epsilon_1 + \ldots + \epsilon_B) } \left( \prod_{1 \leq i \leq B} \frac{ P- k_i}{q} \right)
P^{B - B(Q'_0 - \varepsilon_0)/r_1}
\\
&\ll& P^{B \theta_0} P^{B - B(Q'_0 - \varepsilon_0)/r_1}.
\notag
\end{eqnarray}

Thus we obtain that (\ref{bound in S with linear 0}) is bounded by
\begin{eqnarray}
|S(\boldsymbol{\alpha})| &\ll& q^n \left( \frac{P}{q} \right)^{n - B}  P^{ B \theta_0} P^{B - B(Q'_0 - \varepsilon_0)/r_1}
\notag
\\
&\leq&
q^B P^{n +  B \theta_0  - B(Q'_0 - \varepsilon_0)/r_1}
\notag
\\
&\leq&
P^{n +  2 B \theta_0   - B(Q'_0 - \varepsilon_0)/r_1}.
\notag
\end{eqnarray}
Since we chose $Q$ to satisfy
$$
Q <  B \left(   \frac{Q_0/2 - \varepsilon_0}{r_1} - 2\theta_0 \right),
$$
it follows that we are in alternative $(i)$ as long as $P$ is sufficiently large with respect to
$\mathbf{u}$, $d$, $n$, $r_d, \ldots, r_1$, and $Q$.
\end{proof}

Let $1 < \ell \leq d$ and $r_{\ell} > 0$. We define $g_{\ell}( \mathbf{U}_{\ell} )$
to be the largest real number such that
\begin{equation}
\label{def gd}
z_P(\mathbb{M}_{\ell}) \ll P^{n({\ell}-1) - g_{\ell}( \mathbf{U}_{\ell} ) + \varepsilon}
\end{equation}
holds for each $\varepsilon >0$.
%It was proved in \cite[pp. 280, Corollary]{S} that
%\begin{equation}
%\label{h and g}
%h_{\ell}( \mathbf{U}_{\ell} ) < \frac{\ell!}{  (\log 2)^{\ell} }  \left( g_{\ell}( \mathbf{U}_{\ell}  ) + ({\ell}-1)r_{\ell} (r_{\ell} - 1)  \right).
%\end{equation}
Let
$$
\gamma_{\ell} = \frac{2^{{\ell}-1} ({\ell}-1) r_{\ell}}{ g_{\ell}( \mathbf{U}_{\ell} ) }
$$
when $r_{\ell} >0$ and $g_{\ell}( \mathbf{U}_{\ell} ) > 0$.
We let
$\gamma_{\ell} = 0$ if $r_{\ell} = 0$, and let
$\gamma_{\ell} = + \infty$ if $r_{\ell} > 0$ and $g_{\ell}( \mathbf{U}_{\ell} ) = 0$.
For $\ell$ with $r_{\ell}>0$, we also define
\begin{equation}
\label{def gamma'}
\gamma'_{\ell} = \frac{ 2^{{\ell}-1} }{ g_{\ell}( \mathbf{U}_{\ell} ) } = \frac{ \gamma_{\ell} }{ ({\ell}-1) r_{\ell} }.
\end{equation}

From Lemma \ref{Lemma 15.1 in S}, we obtain the following corollary which is the inhomogeneous polynomials version of \cite[pp.276, Corollary]{S},
and it is obtained by essentially the same proof.
\begin{cor}\cite[pp.276, Corollary]{S}
\label{cor 15.1 in S}
Suppose $\mathbf{u}$ has coefficients in $\mathbb{Z}$.
Let $Q > 0$ and $\varepsilon >0$. Let $2 \leq \ell \leq d$ with $r_{\ell} > 0$. Let $P$ be sufficiently large with respect
to $d$ and $r_d, \ldots, r_1$.
If $\ell =d$, then let $\theta = 0$ and $q=1$.
On the other hand, if $2 \leq \ell < d$, then suppose $0 \leq \theta < 1/4$ and that there is $q \in \mathbb{N}$ with
$$
q \leq P^{\theta} \ \ \text{ and } \ \  \| q \boldsymbol{\alpha}_{j}  \| \leq P^{\theta - j} \ \ (\ell < j \leq d).
$$
Let $S( \boldsymbol{\alpha})$ be the sum associated to $\mathbf{u}$ as in ~(\ref{def of S 1}). %Given $\eta > 0$ such that $\eta + 4 \theta \leq 1$,
Suppose
$$
4 \theta + Q \gamma'_{\ell} < 1.
$$
Then one of the following two alternatives must hold:

$(i)$ $|  S( \boldsymbol{\alpha})  | \leq P^{n-Q}$.

$(ii)$ There exists $n_0 \in \mathbb{N}$ such that
$$
n_0 \ll P^{ Q \gamma_{\ell} + \varepsilon } \text{  and  } \|  n_0 q \boldsymbol{\alpha}_{\ell} \| \ll P^{ -\ell + 4 \theta + Q \gamma_{\ell} + \varepsilon}.
$$
\newline
%When $\ell = d$, the implicit constants
%depend at most on $n,d, r_d, \eta, \varepsilon$ and $\mathbf{U}^{(d)}$.
%Otherwise, when $2 \leq \ell < d$,
The implicit constants
depend at most on $n,d, r_d, \ldots, r_1,\varepsilon$, and $\mathbf{u}$.
\end{cor}

The above corollary does not deal with the case $\ell = 1$, and  we take care of this in the following lemma.
\begin{lem}\cite[Lemma 15.2]{S}
\label{lemma 15.2 in S}
Suppose $\mathbf{u}$ has coefficients in $\mathbb{Z}$, and that
$$
\mathcal{B}_1(\mathbf{u}_1) > 2 r_1 \left( \max \Big{ \{ } \
4 (r_1 + 1) \left( \sum_{j = 2}^d 4^{j- 2} \gamma_j \right)  , \  \frac{1}{4(R+1)} \
\Big{ \} }
\right)^{-1}.
$$
Let $\varepsilon >0$ be sufficiently small. Let $Q > 0$ satisfy
$$
Q (8 r_1 + 8) \left( \sum_{j = 2}^d 4^{j- 2} \gamma_j \right) < 1
\  \   \text{  and  } \  \
\frac{Q}{2 (R + 1)} < 1.
$$
Let $S( \boldsymbol{\alpha})$ be the sum associated to $\mathbf{u}$ as in ~(\ref{def of S 1}).
%Given $\eta > 0$ such that $\eta + 4 \theta \leq 1$,
%Let $\varepsilon_0 > 0$ be sufficiently small. Suppose  $Q > 0$ and $Q_0>0$ satisfy
%$$Q <  ( \mathcal{B}_1(\mathbf{u}_1) - 1) \left(   \frac{Q_0 - \varepsilon_0}{r_1}  - 1 \right) - d.$$
Suppose $P$ is sufficiently large with respect
to $d$, $n$, $r_d, \ldots, r_1$, $\varepsilon$, $Q$, and $\mathbf{u}$.
Then one of the following two alternatives must hold:

$(i)$ $|  S( \boldsymbol{\alpha})  | \leq P^{n-Q}$.

$(ii)$ There exist $n_1, n_2, \ldots , n_d \in \mathbb{N}$ such that
$$
n_{\ell} \ll P^{ Q \gamma_{\ell} + \varepsilon } \ \   \text{  and  }  \  \   \|  n_d \ldots  n_{\ell} \boldsymbol{\alpha}_{\ell} \|
\ll P^{ -\ell + Q \left( \sum_{j = \ell}^d 4^{j- \ell} \gamma_j \right)  + \varepsilon} \ \ (2 \leq \ell \leq d),
$$
$$
n_1 \leq P^{ M_0 Q} \ \ \text{ and }  \  \  \| n_1  \boldsymbol{\alpha}_1  \| \leq P^{- 1 + M_0 Q },
$$
where
$$
M_0 = \max \Big{ \{ } \
8 (r_1 + 1) \left( \sum_{j = 2}^d 4^{j- 2} \gamma_j \right)  , \  \frac{1}{2(R+1)} \ \Big{ \} }.
$$
\newline
The implicit constants depend at most on $n,d, r_d, \ldots, r_1, \varepsilon$, and $\mathbf{u}$.
\end{lem}

\begin{proof} We begin by proceeding as in the proof of \cite[Lemma 15.2]{S}.
Suppose we have
$$
|S(\boldsymbol{\alpha})| > P^{n - Q}.
$$
Let $\varepsilon_d > 0$ be sufficiently small. Since $Q \gamma'_d < 1$, by Corollary \ref{cor 15.1 in S} there exists $n_d \in \mathbb{N}$ with
$$
n_d \ll P^{Q \gamma_d + \varepsilon_d}  \ \ \text{   and  }  \  \   \| n_d  \boldsymbol{\alpha}_d \| \ll  P^{-d + Q \gamma_d + \varepsilon_d}.
$$
Suppose now that $r_{d-1} > 0$. Since $4 Q \gamma_d + Q \gamma'_{d-1} < 1$,
we can apply Corollary \ref{cor 15.1 in S} again with
$\ell = d-1$, $\theta = Q\gamma_d + 2\varepsilon_d$, and $q = n_d$. Note we have by our assumption on $Q$ that $\theta < 1/4$.
Let $\varepsilon_{d-1} > 0$ be sufficiently small. Thus there exists $n_{d-1} \in \mathbb{N}$ with
\begin{equation}
\label{inequalities 123345}
n_{d-1} \ll P^{Q \gamma_{d-1} + \varepsilon_{d-1}}  \ \ \text{   and  }  \  \   \| n_d n_{d-1} \boldsymbol{\alpha}_{d-1} \|
\ll  P^{-(d-1) + 4 Q \gamma_d + 8\varepsilon_d + Q \gamma_{d-1} + \varepsilon_{d-1}}.
\end{equation}
In the case $r_{d-1} = 0$, we have $\gamma_{d-1} = 0$ and obtain (\ref{inequalities 123345})
trivially with $n_{d-1} = 1$. It is clear we can continue in this manner. By repeating the argument, we ultimately obtain
that there exist $n_2, \ldots , n_d \in \mathbb{N}$ such that
$$
n_{\ell} \ll P^{ Q \gamma_{\ell} + \varepsilon } \text{  and  } \|  n_d \ldots  n_{\ell} \boldsymbol{\alpha}_{\ell} \|
\ll P^{ -\ell + Q \left( \sum_{j = \ell}^d 4^{j- \ell} \gamma_j \right)  + \varepsilon} \ \ (2 \leq \ell \leq d).
$$

If $r_1 = 0$, then we are done trivially with $n_1 = 1$. Let $r_1 > 0$.
We now apply Lemma \ref{Lemma 15.1 in S-linear} with
$$
\theta_0 =   \left( \sum_{j = 2}^d 4^{j- 2} \gamma_j \right) Q + d\varepsilon  < 1,
$$
where $\varepsilon > 0$ is sufficiently small,
$$ Q_0/2 =
\max \Big{ \{ } \
4 (r_1 + 1) \left( \sum_{j = 2}^d 4^{j- 2} \gamma_j \right) Q , \  \frac{Q}{4(R+1)} \
\Big{ \} } < \frac12,
$$
and
$$q = (n_d \ldots  n_{2}) \leq P^{\theta_0},$$
where the last inequality holds for $P$ sufficiently large.
Let $\varepsilon_0>0$ be sufficiently small.
With these choices of $\theta_0$ and $Q_0$, we have
$$
2 \theta_0 < (Q_0/2 - \varepsilon_0)/ (2 r_1) < (Q_0/2 - \varepsilon_0)/ r_1.
$$
With our assumption on $\mathcal{B}_1(\mathbf{u}_1)$, it is clear that we have
$$
Q < \mathcal{B}_1(\mathbf{u}_1) \left(  \frac{Q_0/2 - \varepsilon_0}{2 r_1} \right) < \mathcal{B}_1(\mathbf{u}_1) \left( \frac{Q_0/2 - \varepsilon_0}{ r_1} - 2 \theta_0 \right).
$$
Therefore, it follows by Lemma \ref{Lemma 15.1 in S-linear} that there exists $n_1 \in \mathbb{N}$ such that
$$
n_1 \leq P^{Q_0} \ \ \text{ and }  \  \  \|  n_1  \boldsymbol{\alpha}_1  \| \leq P^{Q_0 - 1}.
$$
\end{proof}

%Note we could define $\varepsilon_0$ explicitly in terms of $B_1(\mathbf{u}_1)$ and $r_1$, where
%it makes sure that $Q < \mathcal{B}_1(\mathbf{u}_1) \left(  \frac{Q_0/2 - \varepsilon_0}{2 r_1} \right)$
%as $P$ is supposed to be sufficiently large with respect to $\varepsilon_0$ also otherwise...

We are now in position to prove our main result.
\begin{proof}[Proof of Theorem \ref{prop II in S}]
By the hypotheses, we know that
$$
(8 r_1 + 8) \Delta \Omega \left(  \gamma_2 + 4 \gamma_3 + 4^2 \gamma_4 + \ldots + 4^{d - 2} \gamma_d   \right) < 1,
$$
$$
\frac{\Delta \Omega}{2(R + 1)} < 1,
$$
and
\begin{equation}
\label{Omega bound 1}
\Omega \left(  \gamma_2 + 4 \gamma_3 + 4^2 \gamma_4 + \ldots + 4^{d - 2} \gamma_d   \right) +  \Omega M_0 < 1,
\end{equation}
where
$$
M_0 = \max \Big{ \{ } \
8 (r_1 + 1) \left( \sum_{j = 2}^d 4^{j- 2} \gamma_j \right)  , \  \frac{1}{2(R+1)} \ \Big{ \} }
$$
as in the statement of Lemma \ref{lemma 15.2 in S}.

Let
\begin{equation}
\label{def ep 0}
\varepsilon'_0 = \frac{1}{2 \Omega} \left(   1 - \Omega \left(  \gamma_2 + 4 \gamma_3 + 4^2 \gamma_4 + \ldots + 4^{d - 2} \gamma_d   \right) -  \Omega M_0 \right).
\end{equation}
We apply Lemma \ref{lemma 15.2 in S} with $Q = \Delta \Omega$.
If the alternative $(i)$ of Lemma \ref{lemma 15.2 in S} holds then we are done. Let us suppose we
have the alternative $(ii)$ of Lemma \ref{lemma 15.2 in S}.
Then for $P$ sufficiently large, we have
$$
q := n_d \ldots n_2 n_1 \leq P^{\Delta \Omega \left( \sum_{j = 2}^d 4^{j- 2} \gamma_j \right) + \Delta \Omega M_0 + \Delta \Omega \varepsilon'_0},
$$
and
$$
\| q \boldsymbol{\alpha}_{\ell} \| \leq P^{- \ell + \Delta \Omega \left( \sum_{j = 2}^d 4^{j- 2} \gamma_j \right) + \Delta \Omega M_0
+ \Delta \Omega \varepsilon'_0 } \ \ \ (1 \leq \ell \leq d).
$$
Since
$$
\Omega \left(  \gamma_2 + 4 \gamma_3 + 4^2 \gamma_4 + \ldots + 4^{d - 2} \gamma_d   \right) +  \Omega M_0 + \Omega \varepsilon'_0 < 1,
$$
we obtain our result.
\end{proof}

%We need the following lemma to obtain estimates on the singular integral.
We prove the following lemma which becomes useful in some applications of the Hardy-Littlewood circle method.
The proof is based on that of \cite[Lemma 8.1]{S}.
Let
$$
\mathcal{I}( \mathfrak{B}_0 ,  \boldsymbol{\tau}) = \int_{\mathbf{v} \in \mathfrak{B}_0 } e \left( \sum_{\ell = 1}^d \sum_{r=1}^{r_{\ell}} \tau_{\ell,r}  \cdot U_{\ell,r}(\mathbf{v})  \right) \ \mathbf{d}\mathbf{v}.
$$
\begin{lem}
\cite[Lemma 8.1]{S}
\label{lemma sing int}
Suppose $\mathbf{u}$ has coefficients in $\mathbb{Z}$, and that
$\mathcal{B}_1(\mathbf{u}_1)$ is sufficiently large with respect to $r_d, \ldots, r_1$, and $d$.
Furthermore, suppose $\gamma_2, \ldots, \gamma_d$ are sufficiently small with respect to $r_d, \ldots, r_1$, and $d$.
Then we have
\begin{equation}
\label{(3.9) is S}
\mathcal{I}( \mathfrak{B}_0,  \boldsymbol{\tau}) \ll \min (1 , |\boldsymbol{\tau}|^{-R - 1} ),
\end{equation}
where the implicit constant depends at most on $n$, $d$, $r_d, \ldots, r_1$, and $\mathbf{U}$.
\end{lem}

\begin{proof}
Given $\mathbf{a} = (\mathbf{a}_d, \ldots, \mathbf{a}_1) \in (\mathbb{Z}/q \mathbb{Z})^R$, where
$\mathbf{a}_{\ell} = (a_{\ell,1}, \ldots, a_{\ell,r_{\ell}}) \in (\mathbb{Z}/q \mathbb{Z})^{r_{\ell}}$ $(1 \leq \ell \leq d)$
and $\gcd(\mathbf{a},q) = 1$, let us define
$$
\widetilde{\mathfrak{M}}_{\mathbf{a},q}((R+2)^{-1}) = \{ \boldsymbol{\alpha} \in [0,1)^R :  \max_{1 \leq r \leq r_{\ell}} | q \alpha_{\ell,r}  - a_{\ell,r} | \  \leq  \  P^{(R+2)^{-1}}/P^{\ell} \ \ (1 \leq \ell \leq d)  \},
$$
and let
$$
\widetilde{\mathfrak{M}} = \bigcup_{ q \leq P^{(R+2)^{-1}} } \bigcup_{ \substack{ \mathbf{a} \in (\mathbb{Z}/q \mathbb{Z})^R \\ \gcd(\mathbf{a},q) = 1 }   }  \widetilde{\mathfrak{M}}_{\mathbf{a},q}((R+2)^{-1}).
$$
Note the boxes $\widetilde{\mathfrak{M}}_{\mathbf{a},q}((R+2)^{-1})$ with $q \leq P^{(R+2)^{-1}}$, $\mathbf{a} \in (\mathbb{Z}/q \mathbb{Z})^R$,
and $\gcd(\mathbf{a},q) = 1$ are disjoint when $P$ is sufficiently large.

Suppose $| \boldsymbol{\tau} |> 2$. Let  $P \mathbf{v} = \mathbf{v}'$ so that we have
$$
\mathcal{I}( \mathfrak{B}_0,  \boldsymbol{\tau})  = \frac{1}{P^n} \int_{P \mathfrak{B}_0} e \left( \sum_{\ell = 1}^d \sum_{r=1}^{r_{\ell}} \beta_{\ell,r}  \cdot U_{\ell,r}(\mathbf{v}')  \right) \ \mathbf{d}\mathbf{v}',
$$
where
\begin{equation}
\label{beta in sing int lem}
\beta_{\ell, r} = \frac{\tau_{\ell, r}}{P^{\ell}}  \ \ (1 \leq \ell \leq d, 1 \leq r \leq r_{\ell}).
\end{equation}

Let $P = | \boldsymbol{\tau} |^{R + 2}$, and consider the exponential sum
$$
S(\boldsymbol{\beta}) =
\sum_{ \mathbf{x} \in P \mathfrak{B}_0 \cap \mathbb{Z}^n }  e \left( \sum_{\ell = 1}^d \sum_{r=1}^{r_{\ell}} \beta_{\ell,r}  \cdot U_{\ell,r}(\mathbf{x})  \right).
$$
Then $\boldsymbol{\beta}$ lies on the boundary of the box $\widetilde{\mathfrak{M}}_{\mathbf{0},1}((R+2)^{-1})$.
Thus for $|\boldsymbol{\tau}|$ sufficiently large, $\boldsymbol{\beta}$ lies on the boundary of the set
$\widetilde{\mathfrak{M}}$, which is precisely the set considered in the alternative $(ii)$ of Proposition \ref{prop II in S} with $\Delta = (R+2)^{-1}$.
Consequently, $\boldsymbol{\beta}$ also lies on the boundary of $[0,1)^R \backslash \widetilde{\mathfrak{M}}$.
Since $|S(\boldsymbol{\alpha}) |$ is a continuous function, we obtain via Theorem \ref{prop II in S}
(with $\Omega = R+1$) that
\begin{equation}
\label{ineq in sing int1}
|S(\boldsymbol{\beta}) | \leq P^{n - (R+2)^{-1} \Omega} = P^n |\boldsymbol{\tau}|^{- \Omega} = P^n |\boldsymbol{\tau}|^{-  R - 1}.
%\leq P^n |\boldsymbol{\tau}|^{- R - 1}.
\end{equation}
Note with the hypothesis of this lemma, we have
\begin{eqnarray}
&&\min \Big{ \{ } \  \frac{1}{ 8r_1 + 9} \left( \sum_{j = 2}^d 4^{j- 2} \gamma_j \right)^{-1}, \  \left(  \frac{1}{2 (R + 1)} + \sum_{j = 2}^d 4^{j- 2} \gamma_j  \right)^{-1}  \ \Big{ \} }
\notag
\\
&=& \left(  \frac{1}{2 (R + 1)} + \sum_{j = 2}^d 4^{j- 2}
\gamma_j  \right)^{-1}
\notag
\\
&>& R + 1,
\notag
\end{eqnarray}
which justifies our application of Theorem \ref{prop II in S} with $\Omega = R+1$.

We also have
\begin{eqnarray}
&&S(\boldsymbol{\beta}) - \int_{P \mathfrak{B}_0} e \left( \sum_{\ell = 1}^d \sum_{r=1}^{r_{\ell}} \beta_{\ell,r}  \cdot U_{\ell,r}(\mathbf{v}')  \right) \ \mathbf{d}\mathbf{v}'
\\
&=&
\notag
\sum_{\mathbf{x} \in [0,P)^n }  \int_{x_1}^{x_1+1} \ldots \int_{x_n}^{x_n+1} e \left( \sum_{\ell = 1}^d \sum_{r=1}^{r_{\ell}} \beta_{\ell,r}  \cdot  U_{\ell,r}(\mathbf{x}) \right)   - e \left( \sum_{\ell = 1}^d \sum_{r=1}^{r_{\ell}} \beta_{\ell,r}  \cdot  U_{\ell,r}(\mathbf{v}')  \right)  \ \mathbf{d}\mathbf{v}'
\\
&+& O(P^{n-1})
\notag
\\
&\ll&
\notag
P^n \frac{ |\boldsymbol{\tau}|}{P}  + O(P^{n-1})
\\
&\ll&
\notag
P^{n-1}  |\boldsymbol{\tau}|,
\end{eqnarray}
where we applied the mean value theorem and (\ref{beta in sing int lem}) to obtain the second last inequality.
Therefore, it follows that
$$
S(\boldsymbol{\beta}) = P^n \mathcal{I}( \mathfrak{B}_0,  \boldsymbol{\tau})  + O(P^{n-1}  |\boldsymbol{\tau}|).
$$
It is then easy to deduce from (\ref{ineq in sing int1}) that
$$
\mathcal{I}( \mathfrak{B}_0,  \boldsymbol{\tau}) \ll \min\{ 1,  |\boldsymbol{\tau}|^{- R - 1} \}.
$$
\end{proof}


\begin{thebibliography}{9}

%\bibitem{B}  B. J. Birch,\textit{Forms in many variables}. Proc. Roy. Soc. Ser. A 265 1961/1962, 245--263.

%\bibitem{BGS} J. Bourgain, A. Gamburd and P. Sarnak, \textit{Affine linear sieve, expanders, and sum-product}. Invent. Math. 179 (2010), no. 3, 559--644.

%\bibitem{BHB} T.D. Browning, and D.R. Heath-Brown, \textit{Froms in many variables and differing degrees}. J. Eur. Math. Soc., to appear.

\bibitem{BP} T.D. Browning, and S.M. Prendiville, \textit{Improvements in Birch's theorem on forms in many variables}.
J. Reine Angew. Math., to appear.

%\bibitem{BDLW} J. Br\"{u}dern, R. Dietmann, J. Liu and T. D. Wooley, \textit{A Birch-Goldbach theorem}. Arch. Math. (Basel) 94 (2010), no. 1, 53--58.

\bibitem{CM} B. Cook and {\'A}.  Magyar, \textit{Diophantine equations in the primes}. Invent. Math. {198} (2014), 701--737.

%\bibitem{C} S. Chow, \textit{Roth-Waring-Goldbach}. arXiv:1602.04012.

%\bibitem{D} H. Davenport, \textit{Analytic methods for Diophantine equations and Diopantine inequalities}. Second edition. Cambridge University Press, Cambridge, 2005.

%\bibitem{DRS} W. Duke, Z. Rudnick and P. Sarnak, \textit{Density of integer points on affine homogeneous varieties}. Duke Math. J. 71 (1993), no. 1, 143--179.

%\bibitem{GPY}  D. A. Goldston,  J. Pintz and  C. Y. Y{\i}ld{\i}r{\i}m, \textit{Primes in tuples. I}.  Ann. of Math. (2) 170 (2009), no. 2, 819--862.

%\bibitem{GS}  A. S. Golsefidy and P. Sarnak, \textit{The affine sieve}, J. Amer. Math. Soc. 26 (2013), no. 4, 1085–1105.

%\bibitem{GT1} B. Green and  T. Tao, \textit{The primes contain arbitrarily long arithmetic progressions}. Ann. of Math. (2) 167 (2008), no. 2, 481--547.

% \bibitem{GT} B. Green and  T. Tao, \textit{Linear equations in primes}.  Ann. of Math. (2) 171 (2010), no. 3, 1753--1850.

%\bibitem{GT1} B. Green and  T. Tao, \textit{The distribution of polynomials over finite fields, with applications to the Gowers norms}, Contrib. Discrete Math. 4 (2009), no. 2, 1–36.

%\bibitem{H1} H. A. Helfgott, \textit{Major arcs for Goldbach's problem}. arXiv:1305.2897.

%\bibitem{H2} H. A. Helfgott, \textit{Minor arcs for Goldbach's problem}. arXiv:1205.5252.

%\bibitem{KL} T. Kaufman and S. Lovett, \textit{Worst case to average case reductions for polynomials}, 49th Annual IEEE Symposium on Foundations of Computer Science (2008), 166-175.

%\bibitem{H}  L. K. Hua,  \textit{Additive theory of prime numbers}. Translations of Mathematical Monographs, Vol. 13 American Mathematical Society, Providence, R.I. (1965).

%\bibitem{KW} A. V. Kumchev and T. D. Wooley, \textit{On the Waring-Goldbach problem for eighth and higher powers}. J. London Math. Soc. (2) 93 (2016), no. 3, 811--824.

%\bibitem{L} J. Liu, \textit{Integral points on quadrics with prime coordinates}. Monatsh. Math. 164 (2011), no. 4, 439--465.

%\bibitem{LS}  J. Liu and P. Sarnak, \textit{Integral points on quadrics in three variables whose coordinates have few prime factors}. Israel J. of math., 178 (2010), 393--426.

%\bibitem{L1} Z. Liu, \textit{Small Prime Solutions to Cubic Diophantine Equations}, Canad. Math. Bull. 56(2013), 785-794.

%\bibitem{M}  J. Maynard, \textit{Small gaps between primes}. Ann. of Math. (2) 181 (2015), no. 1, 383--413.

%\bibitem{DamS} D. Schindler, \textit{A variant of Weyl's inequality for systems of forms and applications},	arXiv:1403.7156.

\bibitem{S} W.M. Schmidt, \textit{The density of integer points on homogeneous varieties}.
Acta Math. {154} (1985), no. 3-4, 243--296.

%\bibitem{V} I. M. Vinogradov. \textit{Representation of an odd number as a sum of three primes}. Dokl. Akad. Nauk. SSR, 15:291--294, 1937

\bibitem{XY} S. Y. Xiao and S. Yamagishi, \textit{Zeroes of polynomials in many variables with prime inputs}. arXiv:1512.01258.

%\bibitem{Z}  Y. Zhang, \textit{Bounded gaps between primes}. Ann. of Math. (2) 179 (2014), no. 3, 1121--1174.

\end{thebibliography}
\end{document}